\documentstyle{amsppt}
%
\magnification =\magstep 1
\NoBlackBoxes
\topmatter
\title Hamlet and Pfister forms$^*$\qquad\qquad  (A tragedy in four
acts)\qquad \endtitle

\author J\' an Min\' a\v{c}\endauthor

\abstract  In the mid-1960s A.~Pfister discovered extraordinary,
strongly multiplicative forms which are now called Pfister forms.
From that time on, these forms played a dominant role in the
theory of quadratic forms. One of the key properties of a Pfister
form $q$ is that $q$ extended to a suitable transcendental
extension, has the polynomial $q$ as its similarity factor.
Pfister's original proof used clever matrix calculations.  Here
we show that the desired isometry is induced by the multiplication
of a suitable field element. 

We further consider the surprising possibility that Pfister's
forms were already known by Hamlet, Rosencrantz and Guildenstern,
and that they in fact led to a terrible tragedy which is yet
filled with a haunting beauty and mystery that can still inspire
us to this day. 
  \endabstract


\endtopmatter
\nologo

\document
\hsize=4.8in
\noindent{\bf THE CAST OF CHARACTERS:}
\medskip
\smallskip
\noindent CLAUDIUS\dotfill J\' an Min\' a\v c
\smallskip
\noindent HAMLET\dotfill J\' an Min\' a\v c
\smallskip
\noindent THE QUEEN\dotfill Leslie Hallock
\medskip
\medskip
\noindent  $^*${\it{Performed at Oberwolfach in May, 1992.}}
\vskip.2in
\hoffset=.5in
\hsize=5.8 truein
\parskip=1em
\parindent=0 truein
\item{ACT I\quad} {\bf Hamlet's Terrible Revenge - The Question}
(TRUMPET sounds)
\vskip.2in
\item{THE QUEEN \quad} Something is rotten in the state of Denmark.  How 
weary, stale, flat and unprofitable seem to me all uses of this world!

(walking, throwing his hands up in despair)
\item{HAMLET\quad} Evil is here, I swear it!  In the midst of us.  I 
will avenge it! I will kill neither you nor the King, my dear mother.  
That would be too easy for you!  I will not let you sleep, eat, nor rest.
I will torture you with the mystery of Pfister forms!

\item{THE QUEEN\quad} (falling on her knees, cries out)

Oh!  Have mercy, Hamlet!  Please don't!

\newpage

\item{HAMLET\quad} Frailty, thy name is woman!  Only a solution of this
great puzzle of quadratic forms will free you from my spell!
\item{THE QUEEN\quad} Forms!?  Did you say quadratic forms!?  Hamlet, don't
we have enough forms!?  We already have tax forms, insurance forms, 
government forms, claim forms, request forms; even execution forms - I 
am not filling out anything!  If you want some form in any case, then 
have a uniform!  You can be a soldier!
\item{HAMLET\quad} No, mother.  You will not get rid of me so easily!
Let us waste no time!  Let us share the wisdom of my prophet Horatio.

Assume that $F$ is a field with $1 + 1\neq 0$ .  $K = 
F(X_1,\dots,X_{2n})$, $\varphi$ is an anisotropic quadratic form over $F$,
$\varphi = \langle 1,a_1\rangle \otimes\cdots\otimes \langle 
1,a_n\rangle$.  $\varphi_K$ is the extended form of $F$ to $K$.

As my prophet told me, somewhere in Germany a few centuries later, a man 
named Pfister will discover that 
$\varphi_K\cong\varphi(X_1,\dots,X_{2^n})\varphi_K$.
\item{THE QUEEN\quad} Oh! I am lost! I am utterly lost!  Words!  Words!  
Words!  I don't understand them.
\item{HAMLET\quad} You!  You who married two days after my father's 
burial!  You who married the brother of your husband.  You who doesn't 
have any decency!  YOU LISTEN!
\item{THE QUEEN\quad} You are mad!  Oh, my lord!  (afraid)
\item{HAMLET\quad} Pfister will publish a very beautiful matrix proof, 
and less than 2 years later, Witt with his usual wit and magic, will 
utter a magic word, ``Runde", and the whole proof will collapse into a 
few lines.  It will be a joke.
\item{THE QUEEN\quad} Brevity is the soul of wit!
\item{HAMLET\quad} And also of Pfister.  They always produce elegant and 
short proofs.  But nevertheless, for some time there will be no proof
which uses some underlying algebra with its multiplication and norm-like
map as can be done in the case of $2-, 4-, 8-$, and $16$ dimensional
Pfister forms using quadratic extensions, quaternion and Cayley-Dickson 
algebras respectively.
\item{THE QUEEN\quad} But our good Horatio tells me that Shapiro will make 
some progress.
\item{HAMLET\quad} Yes.  He will.  He will succeed in a quite intriguing 
construction of Pfister forms using the multiplication in Clifford 
algebras.  However, for some time it will not be clear how to use his 
results to obtain multiplicative properties.

\noindent AND THIS IS YOU AND YOUR HUSBAND'S TASK, MY DEAR MOTHER!  YOU 
WILL BE BEWITCHED, RESTLESS, WORRIED ... UNTIL YOU FIND THIS PROOF!  
THIS IS MY REVENGE!

\newpage

\vskip.2in
\item{\quad} TRUMPET\qquad  FLASH CARD
\vskip.2in
\item{ACT II\quad} {\bf A Quest for an Algebra and a Marriage in Danger}
\item{THE QUEEN\quad} Claudius, you are a king; you have power, money, gold; 
everything!  You \underbar{must} find an answer!  What would be the 
right generalization of a quadratic field, a quaternion algebra, or a 
Cayley-Dickson algebra?
\item{THE KING\quad} (sadly) My dear Queen, unfortunately there is no king's 
road to mathematics.  Horatio even told me that Hurwitz and other 
mathematicians will prove that no 16 dimensional algebra will work.  We 
are lost.

\item{THE QUEEN\quad} I should have never married you.  We should have never 
made a union.  We cannot even prove a theorem together!
\item{THE KING\quad} We shall, we will, we must!

Look, a 1-fold Pfister form $\langle 1,a\rangle$ is a norm form from the 
single quadratic extension $L = F(\sqrt{-a})$, so if there is any justice
in the world, then an $n$-fold Pfister form must be the norm from the 
multiquadratic extension $L = F(\sqrt{-a_1},\dots,\sqrt{-a_n})$.
\item{THE QUEEN\quad} By my lord, don't you know that there is no justice in 
the world?  It will not work!
\item{THE KING\quad} (whispers)  I found the following interesting scrap of 
paper in Hamlet's chamber.  Poor Hamlet, he doesn't know that I have all 
of the keys to this castle (still looking at various scraps of paper).
A love letter to Ophelia - the idiot!  (He throws the letter away.)

Look, look!  Here is something of interest.

$L = K\left(\sqrt{-a_1}, \sqrt{-a_2\dsize\frac{\widehat\psi_1}{\psi_1}},\cdots,
\sqrt{-a_n\dsize\frac{\widehat\psi_{n-1}}{\psi_{n-1}}}\right)$ where

$\psi_1 = X^2_1 + a_1X^2_2,\;\; \widehat\psi_1 = X^2_3 + a_1X^2_4$,

$\psi_2 = X^2_1 + a_1X^2_2 + a_2X^2_3 + a_1a_2X^2_4,\;\; \widehat\psi_2 =
X^2_5 + a_1X^2_6 + a_2X^2_7 + a_1a_2X^2_8$, etc.  Observe that 
$L\cong\dsize\bigotimes^n_{i=1} L_i$ as $K$-vector spaces, where $L_1 =
K(\sqrt{-a_1})$,

$L_i =
K\left(\sqrt{-a_i\dsize\frac{\widehat\psi_{i-1}}{\psi_{i-1}}}\right)$ for 
$i\geq 2$.  $\ell_1\otimes\cdots\otimes \ell_n\to \ell_1\cdots\ell_n$ as
any fool will tell you.  This must be our algebra!

\newpage

\item{THE QUEEN\quad} You are a thief - a scoundrel!  How could you do it?
How could you steal!?
\item{THE KING\quad} Thus conscience does make cowards of us all.  But I 
have no conscience; so I am not a coward; therefore I steal.
Alright, look!
\input table.sty
\table{ll}
$\kern-.5in L = K(\sqrt{-a}),\; \omega = N_{L/K}$-1-fold Pfister
 &\quad  $L  = K\left(\sqrt{-a_1},
   \sqrt{-a_2\dsize\frac{\widehat\psi_1}{\psi_1}},\cdots,
   \sqrt{-a_n\dsize\frac{\widehat\psi_{n-1}}{\psi_{n-1}}}\right)$ \\ \\
 $\kern-.5in \Theta = X_1 + X_2\sqrt{-a}$ & \quad $\omega =$ ? \qquad   $\Theta =$ ? \\ \\
$\kern-.5in \omega(Z_1 + Z_2\sqrt{-a})$ is \\
 $\kern-.5in Z^2_1 + Z^2_2a = N_{L/K}(Z_1 + Z_2\sqrt{-a})$ &\quad
   $\omega(\Theta\cdot\ell) = \psi(X_1,\cdots, X_{2^n})\omega(\ell)$ \\ \\
  $\kern-.5in \omega(\Theta\ell)  = \omega(\Theta)\omega(\ell)$   & \quad
$\psi$ is just $\varphi$ evaluated at $X_1,\cdots,X_{2^n}, \psi\in K$. \\
$\kern-.5in \hskip.4in = \psi(X_1, X_2)\omega(\ell).$
\endtable

We just have to find a good way of describing $\omega$ on $L$ and also a 
good element $\Theta\in L$.
\item{THE QUEEN\quad}  You are naive, my lord.  This cannot work.  I am sure 
that you need a noncommutative nonassociative algebra; something 
horrible, something beyond our wildest imaginings.  You need more matter 
with less art.
\item{THE KING\quad} How dare you call me naive!  I steal, I murder, I have
power -- I am the king!  I cannot be naive.  I order $\Theta$ to exist!
\item{THE QUEEN\quad} (ironically)  Well, if you are so sure, find it!
\item{KING\quad} (laughing)  I have it!  (ha,ha)  Polonius stole it from
Horatio's room.  Here it is:

$$\Theta = (X_1 + X_2\sqrt{-a_1})\left( 1 +
\sqrt{-a_2\dsize\frac{\widehat\psi_1}{\psi_1}}\right)\cdots \left( 1 + 
\sqrt{-a_n\dsize\frac{\widehat\psi_{n-1}}{\psi_{n-1}}}\right).$$
\item{THE QUEEN\quad} I hate you passionately!
\vskip.2in
\item{} TRUMPET \qquad FLASH CARD
\vskip.2in

\newpage

\item{ACT III\quad} {\bf The Proof}
\item{THE QUEEN\quad} (walking) To prove or not to prove, that is the 
question.
\item{CLAUDIUS\quad} To prove! To prove! But what to prove?
\item{THE QUEEN\quad} You proved that you are a scoundrel and a common 
thief!  That is all that you have achieved.
\item{CLAUDIUS\quad} I know all of this.  But I still love you, still 
adore you.  I need your help.  We have to show that
$$\omega(\Theta\cdot\ell) = \psi(X_1,\dots,X_{2^n})\omega(\ell).$$
Rosencrantz and Guildenstern found out how $\Theta$ should be defined.
$\omega = \dsize\bigotimes^n_{i=1} \omega_i$, where $\omega_i$ is the norm
form $N_{L_{i/ K}}$ on $L_i$.  You remember that $L \cong
\dsize\bigotimes^n_{i=1} L_i$ as $K$-vector spaces.
\item{THE QUEEN\quad} Poor Rosencrantz and Guildenstern!  I heard that they 
were executed in England.
\item{CLAUDIUS\quad} There is a divinity that shapes our ends, rough-hew 
them how we will.  In any case, before their death they managed to tell 
me the proof.  Look!  It works beautifully:

Define $B$ to be a symmetric bilinear form associated with $\omega$ and 
define $\widetilde\omega$ by the equation $\widetilde\omega(\ell) = \omega(\Theta\ell)$.  Let
$B$ be also the symmetric bilinear form associated with $\omega$.  Then 
we have for all $h_1,\dots,h_n$ and $\ell_1,\dots,\ell_n\in L$
$$\align
    \widetilde B(h_1\otimes\cdots\otimes h_n, \ell_1\otimes\cdots\otimes
    \ell_n)
    &=  \widetilde B(h_1\cdots h_n, \ell_1\cdots \ell_n) \\
    &=  B(\Theta h_1\cdots h_n, \Theta\ell_1\cdots \ell_n) \\
    &=  B\left(\left(\dsize\prod^n_{i=1}\Theta_ih_i\right),
    \prod^n_{i=1}\left(\Theta_i\ell_i\right)\right) \\
    &= \dsize\prod^n_{i=1} B_i(\Theta_i h_i, \Theta_i\ell_i) \\
    &= \dsize\prod^n_{i=1} N_{L_i/ K}(\Theta_i)\dsize\prod^n_{i=1}
    B_i(h_i,\ell_i) \\
    &= \psi(X_1,\dots,X_{2^n})B(h_1\otimes \cdots\otimes h_n, 
    \ell_1\otimes\cdots\otimes\ell_n). \endalign $$
Observe

$$\prod^n_{i=1} N_{L_i/ K}(\Theta_i) = \psi(X_1,\dots,X_{2^n}).$$

\newpage

Indeed $N_{L_1/ K}(\Theta_1) = X^2_1 + a_1X^2_2 = \psi_1(X_1,X_2)$.
$$\align N_{L_1/ K}(\Theta_1) N_{L_2/ K}(\Theta_2) & = (X^2_1 +
a_1X^2_2)\left(1 + a_2\dsize\frac{x^2_3 + a_1X^2_4}{X^2_1 + 
a_1X^2_2}\right) \\
 & = \psi_1(X_1,X_2,X_3,X_4) \endalign$$
etc.

\item{THE KING\quad} Finally from the linearity reason we conclude
$$\widetilde B = \psi(X_1,\dots,X_{2^n})B.$$
\item{THE QUEEN\quad} What a pity Rosencrantz and Guildenstern were executed.
They were gifted!
\item{THE KING\quad} Well, this is the reason they were executed.  No decent
monarchy can support creative people.
\vskip.2in
\item{} TRUMPET \qquad FLASH CARD
\vskip.2in
\item{ACT IV\quad} {\bf Pfister Matrices and the Death of a King}
\item{THE QUEEN\quad} So we found that a suitable multiplication in some 
field does the trick.  What is it good for?
\item{THE KING\quad} I guess, if you have gained such good insight as we 
have obtained, you can try to generalize; for example, why not replace 
our field extensions $L_i/K$ by suitable cubic extensions; perhaps one 
can discover cubic Pfister forms, or find nice subspaces of $L/K$, or
study the spaces of similarities; perhaps one can make further progress 
with Shapiro's conjecture.  In particular, perhaps one can even solve 
his ``Pfister-factor" conjecture.
\item{THE QUEEN\quad} Why would one like to solve a conjecture?
\item{THE KING\quad} If it will be solved, it will be discovered now; if it 
won't be discovered, it won't be now; if it isn't now, yet it will be 
solved; readiness is all.
\item{THE QUEEN\quad} What will happen if we can solve it?
\item{THE KING\quad} The cat will mew, and the dog will have his day.
\item{THE QUEEN\quad} I wonder whether poor Rosencrantz and Guildenstern's
proof has something in common with Pfister and Witt's proofs.

\newpage

\item{THE KING\quad} Yes.  They told me that their proof essentially gives a 
basis-free variation of Pfister's proof; it shows that Pfister matrices 
correspond to field multiplication.  When compared with Witt's proof, 
Rosencrantz and Guildenstern's proof not only shows that $\psi_K$ and 
$\psi(X_1,\dots,X_{2^n})\varphi_K$ are isometric, but also constructs 
isometry.
\item{THE QUEEN\quad} Rosencrantz and Guildenstern formed such a cute 
isometric pair!  One did exactly the same as the other.  They were 
indistinguishable.  No wonder they found a good isometry!

(Queen sips from wine cup)
\item{THE KING\quad} Don't drink that wine!  That wine is for Hamlet.
\item{THE QUEEN\quad} O.K., O.K.  I thought that there was enough wine in 
our kingdom.
\item{THE KING\quad} Not any more!  We have to be careful.  The economy is 
very bad.  We have to be as stingy as possible.

In any case, we can pick a good orthogonal ordered basis $B$ of our 
field $L$ over $K$ and identify

\table{ll}
 $L\longleftrightarrow$   &\quad $\left\{{\left(\array{c} Y_1 \\ \vdots \\
 Y_{2^n}\endarray\right)}, Y_i\in K\right\} = V$ \\ \\
 $e = \Sigma Y_ib_i$
 &\quad ${\left(\array{c} Y_1 \\ \vdots \\ Y_{2^n}\endarray\right)}$ \\ 
 \\
 $\omega$ &\qquad $\varphi$ \\ \\
 $\widetilde\Theta \colon L\to L$ &\quad $\widetilde{\widetilde \Theta} 
 \colon V\to V$ \\ \\
$e\to\Theta\ell$ & \quad
  ${\left(\array{c} Y_1 \\ \vdots \\ Y_{2^n}\endarray\right)} \to
{T \left(\array{c} Y_1 \\ \vdots \\ Y_{2^n}\endarray\right)}$ \\
\endtable

\newpage

\table{ll}
& \quad where
$T$ is a matrix of $\widetilde{\widetilde\Theta}$ written with respect to $B$ \\ \\
$\omega(\Theta\ell) = \omega(\Theta)\omega(\ell)$ &\quad $\varphi{\left(T
 \left(\array{c} Y_1 \\ \vdots \\ Y_{2^n}\endarray\right)\right)} =
 \varphi(X_1,\dots,X_{2^n})\varphi(Y_1,\dots, Y_{2^n})$
 \endtable

which shows beautifully that $\varphi$ is strictly multiplicative in the 
Pfister sense.

(The Queen tries to secretly drink the same wine as before.  The King 
observes it.  He rushes to the queen, takes the wine out of her hand.  
He is furious.)
\item{THE KING\quad} I TOLD YOU NOT TO DRINK THAT WINE!!  We can't afford 
it!  Would you like to ruin our whole kingdom?  You will become an 
alcoholic.
\item{THE QUEEN\quad} You horrible man!  (to audience) What a cheap man my 
husband is!  All that I wanted was to sweeten my tongue!  (cynically) 
Lend me this wine, my Lord.  I will return it to you.
\item{THE KING\quad}  Neither a borrower nor a lender be.  NO!

Listen my dear.  Don't drink it.  Rather, pay attention to this
beautiful example:

$\varphi = \langle 1,a_1\rangle\otimes\langle 1,a_2\rangle$.  Then $K = 
F(X_1,X_2,X_3,X_4)$.

\table{ll}
$L = F\left(\sqrt{-a_1},\sqrt{-a_2\dsize\frac{X^2_3 + a_1X^2_4}{X^2_1 +
a_1X^2_2}}\right)$ & \quad
$T ={\left(\array{rrrr}
        X_1 & -aX_2 & S & \\
        X_2 &   X_1 &   & \\
        X_3 & -a_1X_4 & X_1 & -a_1X_2 \\
        X_4 & X_3   & X_2 & X_1\endarray\right)}$  \endtable

$B = \left\{1,\sqrt{-a_1},\sqrt{-a_2\dsize\frac{\hat\psi_1}{\psi_1}}\;\dsize\frac{X_1 +
X_2\sqrt{-a_1}}{X_3 + X_4\sqrt{-a_1}},
\sqrt{-a_1}\sqrt{-a_2\dsize\frac{\hat\psi_1}{\psi_1}}\;\dsize\frac{X_1 +
X_2\sqrt{-a_1}}{X_3 + X_4 \sqrt{-a_1}}\right\}$

where

$$\kern-.5in S = \frac{-a_2}{X^2_1 + a_1X^2_2}
\left(\array{l|l}
 X^2_1X_3 + 2a_1X_1X_2X_4 - a_1X^2_2X_3
 & -2a_1X_1X_2X_3-a^2_1X^2_2X_4 +
 a_1X^2_1X_4 \\ \\ \hline \\
 2X_1X_2X_3-X^2_1X_4 + a_1X^2_2X_4 &
 X^2_1X_3 - a_1X^2_2X_3 +
 2a_1X_1X_2X_4\endarray\right)$$

\newpage

(shows a large paper with $S$ written out)
 
\item{THE QUEEN\quad} Enough mathematics!  I want to be drunk!  (She drains 
the glass, falls to the floor, crying out in pain)  Oh, OH!! (etc.)
\item{THE KING\quad} Gertrude!  What have you done!?  There was poison in 
that wine.

(He puts a sword near Gertrude - leans over her)
\item{THE QUEEN\quad} WHAT!?  (furiously - with her last dying strength)  You
wanted to kill my son?!  You murderer!  You scoundrel!  You miserable 
man!

(The Queen grasps the sword and plunges it into the King - blood flows from the 
King's heart)
\item{THE KING\quad} (Tries to stay on his feet - he staggers)  We shall 
never publish our solution!  PUBLISH OR PERISH!  OH!  OH!

(The King falls to the ground, suffering.  They both are dead)

\medskip

\noindent {\bf ACKNOWLEDGEMENTS}

\noindent Hamlet and the Queen gratefully
acknowledge the continuous interest of Susanne Pumpl\"un in this play,
and its consequences. We also thank her for her correction of basis $B$
on page 8. 

\noindent {\bf Remark} (July 2009). Karim Johannes Becher proved
the Pfister factor conjecture in {\it{Invent. Math.}} {\bf 173} (2008),
1-6.

\noindent Susanne Pumpl\"un was able to extend the ideas in this
play to higher dimensional forms (see her forthcoming papers).

\noindent This play was performed in May 1992 at Oberwolfach with
Albrecht Pfister present! 

\enddocument

\end